\documentclass{article}
\usepackage{amssymb}
\usepackage{graphicx}
\usepackage{amsmath}

\newtheorem{theorem}{Théorème}[section]

\newtheorem{conclusion}[theorem]{Conclusion}

\newtheorem{corollary}[theorem]{Corollaire}

\newtheorem{definition}[theorem]{Définition}
\newtheorem{example}[theorem]{Exemple}

\newtheorem{notation}[theorem]{Notation}

\newtheorem{proposition}[theorem]{Proposition}
\newtheorem{remark}[theorem]{Remarque}

\newenvironment{proof}[1][Preuve]{\textbf{#1.} }{\ \rule{0.5em}{0.5em}}
\input{tcilatex}

\begin{document}

\title{{\LARGE Extension d'un feuilletage de Lie minimal }\\
{\LARGE \ d'une vari\'{e}t\'{e} compacte}}
\author{\textit{Cyrille Dadi }$^{(1)}$\textit{\ et Hassimiou Diallo }$^{(2)}$
\\
%EndAName
\textit{Laboratoire de Math\'{e}matiques Fondamentales,} \\
\textit{Ecole Normale Sup\'{e}rieure, Universit\'{e} de Cocody}\\
\textit{08 BP 10 ABIDJAN 08}\\
\textit{email}$^{(1)}$\textit{: cyrdadi@yahoo.fr, email}$^{(2)}$\textit{:
diallomh@yahoo.fr}}
\maketitle

\begin{abstract}
\textit{The purpose of this paper is to show that any extension of a minimal
\ Lie foliation on a compact manifold is a transversely Riemannian }$\frac{%
\mathcal{G}}{\mathcal{H}}-$\textit{foliation with trivial normal bundle. }

\textit{This result permits to classify the extensions of a minimal Lie
foliation on a compact manifold from Lie subgroups of its Lie group.}
\end{abstract}

Dans tout ce qui suit les objets sont $C^{\infty },$ la vari\'{e}t\'{e} de
base $\mathcal{M}$ et les feuilletages consid\'{e}r\'{e}s sont, si
necessaire, pris orientables, et le groupe de Lie $G$ des feuilletages de
Lie sera suppos\'{e} connexe et simplement connexe.

\section{Introduction}

On sait qu'un feuilletage de Lie d'une vari\'{e}t\'{e} compacte et connexe
est , entre autre, un feuilletage transversalement - riemannien , parrall%
\'{e}lisable, homog\`{e}ne-, \`{a} fibr\'{e} normal trivial , d\'{e}fini par
une une $1-forme$ \`{a} valeurs dans une alg\`{e}bre de Lie. \newline
Il s'agit ici , \'{e}tant donn\'{e} un $G-$feuilletage de Lie \textit{minimal%
} $\mathcal{(M}$, $\mathcal{F)}$ d'une vari\'{e}t\'{e} compacte et connexe ,
de voir comment \'{e}voluent ces propri\'{e}t\'{e}s \`{a} travers des
feuilletages $\mathcal{F}^{\prime }$\textit{extensions} de $\mathcal{F}$ ( $%
i.e.$des feuilletages qui sont tels que \newline
$T\mathcal{F}$ $\varsubsetneq T\mathcal{F}^{\prime }$). A d\'{e}faut de pr%
\'{e}server ces propri\'{e}t\'{e}s sus-cit\'{e}es, peut-on expliciter une
condition n\'{e}cessaire et suffisante pour qu'une extension admette comme
structure transverse l'une des structures transverses de $\mathcal{F}$?

L'\'{e}tude fait appara\^{\i}tre clairement que le\textit{\ groupe de Lie }$%
G $ contient toutes les\textit{\ informations} concernant l'existence et la
nature d'une telle extension.

De fa\c{c}on pr\'{e}cise , on montre que

\textit{- il y a une correpondance biunivoque entre les sous-alg\`{e}bres de
Lie de }$\mathcal{G=}Lie(G)($\textit{\ ou si l'on pr\'{e}f\`{e}re entre les
sous -groupes de Lie connexes de }$G)$\textit{\ et les extensions de }$%
\mathcal{F}$\textit{,}

\textit{- une extension de }$\mathcal{F}$\textit{\ est un }$\frac{\mathcal{G}%
}{\mathcal{H}}-$\textit{feuilletage ( voir d\'{e}finitions) transversalement
riemannien , \`{a} fibr\'{e} normal trivial , d\'{e}fini par une }$1-$%
\textit{\ forme vectorielle.}

Ce qui permet d'obtenir une caract\'{e}risation et par suite une
classification des extensions de $\mathcal{F}$.

\textit{Une extension }$\mathcal{F}$\textit{' d'un feuilletage de Lie
minimal }$\mathcal{(M}$, $\mathcal{F)}$\textit{\ d'une vari\'{e}t\'{e}
compacte et connexe est transversalement homog\`{e}ne (resp. de Lie) si et
seulement si le sous-groupe de Lie connexe correspondant est ferm\'{e}
(resp. normal).}

Il en r\'{e}sulte que

\textit{- toute extension d'un feuilletage de Lie (resp.d'un feuilletage lin%
\'{e}aire)}\newline
\textit{\ minimal du tore est un feuilletage de Lie (resp. un feuilletage lin%
\'{e}aire) et}

\textit{- si un feuilletage de Lie d'une vari\'{e}t\'{e} compacte est dense
dans une de ses extensions alors cette extension est un }$\frac{\mathcal{G}}{%
\mathcal{H}}-$\textit{\ feuilletage transversalement riemannien \`{a} fibr%
\'{e} normal trivial.}

\section{G\'{e}n\'{e}ralit\'{e}s}

\begin{notation}
Dans ce qui suit $\mathcal{G}$ \ est une alg\`{e}bre de Lie de dimension q,
de groupe de Lie connexe et simplement connexe $G,$ $\mathcal{H}$ une
sous-alg\`{e}bre de Lie de $\mathcal{G}$ \ de codimension q', ($%
e_{1},...,e_{q})$ une base de l'alg\`{e}bre de Lie $\mathcal{G}$ \ telle que
($e_{q^{\prime }+1},...,e_{q})$ soit une base de $\mathcal{H}$; on pose [$%
e_{i},e_{j}]=\sum\limits_{k=i}^{q}K_{ij}^{k}e_{k}$ $,$ les $K_{ij}^{k}$ \'{e}%
tant les constantes de structure de $\mathcal{G}$. Ainsi si $\omega $ est
une 1-$forme$ sur une vari\'{e}t\'{e} $\mathcal{M}$ \ \`{a} valeurs dans $%
\mathcal{G}$, relativement \`{a} cette base ,\ on a $\ \omega
=\sum\limits_{i=1}^{q}\omega ^{i}\otimes e_{i}$ qu'on note encore $\omega
=(\omega ^{1},...,\omega ^{q});$ par exemple si $\theta $ est la 1-$forme$
canonique de G, on \'{e}crira $\theta =\sum\limits_{i=1}^{q}\theta
^{i}\otimes e_{i}$ ou $\theta =(\theta ^{1},...,\theta ^{q}).$
\end{notation}

Pr\'{e}cisons pour la suite qu'un sous-groupe de Lie \'{e}tant vu comme un
sous-groupe\textit{\ immerg\'{e}} d'un groupe de Lie, n'est alors ni
necessairement \textit{plong\'{e},} ni necessairement \textit{ferm\'{e}}.

\begin{definition}
Une extension d'un feuilletage ( $\mathcal{M}$, $\mathcal{F}$) de
codimension q est un feuilletage ($\mathcal{M}$,$\mathcal{F}$' ) de
codimension q'tel que 0\TEXTsymbol{<} q'\TEXTsymbol{<}q et T$\mathcal{F}$ $%
\subset $T$\mathcal{F}$' \newline
( on notera $\mathcal{F}\subset $ $\mathcal{F}$').

Une extension d'un feuilletage sera dite de Lie, resp.tranversalement homog%
\`{e}ne, resp. lin\'{e}aire si cette extension est dans chacun des cas un
feuilletage de ce type.
\end{definition}

La notion de $\frac{\mathcal{G}}{\mathcal{H}}-$feuilletage , que nous allons
d\'{e}finir maintenant, a \'{e}t\'{e} introduite par El Kacimi dans \cite%
{EGN}.\newpage

\begin{definition}
Avec les notations ci-dessus , soient $\mathcal{G}$ une alg\`{e}bre de Lie , 
$\mathcal{H}$ une sous-alg\`{e}bre de Lie de $\mathcal{G}$ et $\omega
=(\omega ^{1},...,\omega ^{q})$ \ une $\mathit{1-}forme$ sur une vari\'{e}t%
\'{e} connexe $\mathcal{M}$ \`{a} valeurs dans $\mathcal{G}$. Supposons que $%
\omega $ v\'{e}rifie la condition de Mauer-Cartan $d\omega +\frac{1}{2}%
[\omega ,\omega ]=0,i.e.$%
\begin{equation*}
d\omega ^{k}=-\frac{1}{2}\sum\limits_{i,j=1}^{q}K_{ij}^{k}\omega ^{i}\wedge
\omega ^{j}\text{ \ }(\ast )
\end{equation*}%
\ et $\omega ^{1},...,\omega ^{q^{\prime }}$ sont lin\'{e}airement ind\'{e}%
pendantes en tout point de $\mathcal{M}$. Alors le syst\`{e}me diff\'{e}%
rentiel $\omega ^{1}=...=\omega ^{q}$ $=0$ est int\'{e}grable et d\'{e}finit
un feuilletage $\mathcal{F}$ de codimension $q^{\prime }$ qu'on appellera un 
$\frac{\mathcal{G}}{\mathcal{H}}-$feuilletage d\'{e}fini \ par la 1-$forme$ $%
\omega .$

Si la 1-$forme$ $\omega $ est la 1-$forme$ de F\'{e}dida d\'{e}finissant un
feuilletage de Lie $\mathcal{F}_{\omega },$ on dira que $\mathcal{F}$ est le 
$\frac{\mathcal{G}}{\mathcal{H}}-$feuilletage associ\'{e} au feuilletage de
Lie $\mathcal{F}_{\omega }.$
\end{definition}

\begin{example}
Si $M=G,$ alors $\theta =(\theta ^{1},...,\theta ^{q})$ d\'{e}finit un $%
\frac{\mathcal{G}}{\mathcal{H}}-$feuilletage $\mathcal{F}_{G,H}$ dont les
feuilles sont les classes \`{a} gauche de $H.$
\end{example}

\bigskip On notera que si le sous-groupe $H$ est ferm\'{e}, ce $\frac{%
\mathcal{G}}{\mathcal{H}}-$feuilletage n'est riemannien que si et seulement
si , les actions \`{a} droite de $H$ sont des isom\'{e}tries pour la
structure m\'{e}trique invariante \`{a} gauche de groupe de Lie de $G$.

Ainisi, pour $G=\mathcal{A}ff(\mathbb{R}^{q})=GL(\mathbb{R}^{q})\ltimes 
\mathbb{R}^{q},$ les seules actions \`{a} droite invariantes \'{e}tant les 
\'{e}l\'{e}ments du groupe orthogonal $O(q,\mathbb{R)}$, alors le
feuilletage $\mathcal{F}_{\mathcal{A}ff(\mathbb{R}^{q}),GL(\mathbb{R}^{q})}$
n'est riemannien pour aucune valeur de $\ q\geq 1$ puisque $O(q,\mathbb{%
R)\varsubsetneq }GL(\mathbb{R}^{q}).$

La proposition qui suit , d'apr\`{e}s \cite{EGN}, s'\'{e}tablit comme pour
les feuilletages de Lie ou les feuilletages transversalement homog\`{e}nes(%
\cite{BLU}, \cite{FED}).

\begin{proposition}
\label{prop.}Soit $\mathcal{F}$ un $\frac{\mathcal{G}}{\mathcal{H}}-$%
feuilletage sur une vari\'{e}t\'{e} $\mathcal{M}$ d\'{e}finie par une 1-$%
forme$ $\omega $ et soit $\widetilde{\mathcal{F}}$=$p^{\ast }\mathcal{F}$ le
feuilletage relev\'{e} de $\mathcal{F}$ sur le rev\^{e}tement universel $%
\widetilde{\mathcal{M}}$ \ de $\mathcal{M}$ . Alors, il existe une
application $\mathcal{D}$: $\widetilde{\mathcal{M}}$ $\mathcal{\rightarrow }%
G $ et \ une repr\'{e}sentation \newline
$\rho :\pi _{1}(\mathcal{M}$ )$\mathcal{\rightarrow }G$ telles que

(i) $\mathcal{D}$ est $\pi _{1}(\mathcal{M}$ )-\'{e}quivariant, $i.e.%
\mathcal{D(\gamma }$.$\widetilde{x})=\rho (\gamma )$.$\mathcal{D(}\widetilde{%
x})$ pour tous $\widetilde{x}\in \widetilde{\mathcal{M}}$ \ et $\gamma \in
\pi _{1}(\mathcal{M}$ ), et

(ii) $\ p^{\ast }\omega =\mathcal{D}^{\ast }\theta $, $i.e.$ $\widetilde{%
\mathcal{F}}$=$\mathcal{D}^{\ast }\mathcal{F}_{G,H}.$
\end{proposition}

On dira que $\mathcal{D}$ est une \textit{application d\'{e}veloppante} sur $%
\widetilde{\mathcal{M}}$ du $\frac{\mathcal{G}}{\mathcal{H}}-$feuilletage $%
\mathcal{F}$.

\begin{remark}
\label{Rem1}R\'{e}ciproquement si l'on se donne une r\'{e}pr\'{e}sentation%
\newline
$\rho :\pi _{1}(\mathcal{M}$ )$\mathcal{\rightarrow }G$ et une submersion $%
\pi _{1}(\mathcal{M}$ )-\'{e}quivariant $\mathcal{D}$ de $\widetilde{%
\mathcal{M}}$ sur $G.$ Alors le feuilletage $\widetilde{\mathcal{F}}$=$%
\mathcal{D}^{\ast }\mathcal{F}_{G,H}$ passe au quotient et d\'{e}finit sur $%
\mathcal{M}$ un \ feuilletage $\mathcal{F}$ qui est un $\frac{\mathcal{G}}{%
\mathcal{H}}-$feuilletage.

En particulier si $\mathcal{D}$ est une d\'{e}veloppante de F\'{e}dida d'un
feuilletage de Lie $\mathcal{F}_{\mathcal{D}}$, alors $\mathcal{F}$ est un $%
\frac{\mathcal{G}}{\mathcal{H}}-$feuilletage extension de $\mathcal{F}_{%
\mathcal{D}}$.

\ \ En plus si le sous-groupe $H$ est ferm\'{e}, cette extension est un
feuilletage riemannien si et seulement si le feuilletage $\mathcal{F}_{G,H}$
est riemannien.
\end{remark}

\begin{example}
\label{Exemple}
\end{example}

Avec les notations pr\'{e}c\'{e}dentes, si $\mathcal{F}_{\mathcal{D}}$ est
l'un des deux fllots "propres" du tore hyperbolique $\mathbb{T}_{A}^{3}$ $%
\cite{CAR},$ le diagramme suivant

\begin{equation*}
\begin{array}{ccccc}
\widetilde{\mathbb{T}_{A}^{\text{ }3}} & \overset{}{\longrightarrow } & 
\mathcal{A}ff^{+}(\mathbb{R})=\mathbb{R}_{+}^{\ast }\ltimes \mathbb{R} & 
\overset{pr_{2}}{\longrightarrow } & \mathbb{R}=\frac{\mathcal{A}ff^{+}(%
\mathbb{R})}{\mathbb{R}_{+}^{\ast }} \\ 
\downarrow &  &  &  &  \\ 
\mathbb{T}_{A}^{3} &  &  &  & 
\end{array}%
\end{equation*}

montre que le $\frac{\mathcal{A}ff^{+}(\mathbb{R})}{\mathbb{R}_{+}^{\ast }}%
-feuilletage$ $\mathcal{F}$ - extension de ce flot de Lie (qui est \textit{%
non minimal} $\cite{CAR}$) est un feuilletage non riemannien.

Rappelons, pour terminer cette partie, qu'un groupe $G$ est dit \textit{%
virtuellement r\'{e}soluble} s'il contient un sous-groupe r\'{e}soluble
d'indice fini ( c'est le cas des groupes r\'{e}solubles, nilpotents et ab%
\'{e}liens). Et le th\'{e}or\`{e}me suivant d\^{u} \`{a} A.Haefliger \cite%
{HAE} nous sera utile pour la suite.

\begin{theorem}
\label{Hae}Un feuilletage riemannien \`{a} feuilles denses sur une vari\'{e}t%
\'{e} riemannienne compl\`{e}te $\mathcal{M}$ \`{a} groupe fondamental $\pi
_{1}(\mathcal{M}$ ) virtuellement r\'{e}soluble est transversalement homog%
\`{e}ne.
\end{theorem}

\section{Extension d'un feuilletage de Lie minimal}

\bigskip La propri\'{e}t\'{e} suivante montre la rigidit\'{e} des extensions
d'un feuilletage de Lie minimal d'une vari\'{e}t\'{e} compacte.

\begin{proposition}
\label{pro.2} \ Si $\mathcal{F}^{\prime }$est une extension d'un feuilletage
de Lie minimal \newline
( $\mathcal{M}$, $\mathcal{F}$) d'une vari\'{e}t\'{e} compacte connexe,
alors tout champ global $\mathcal{F-}$feuillet\'{e} transverse, tangent \`{a}
$\mathcal{F}^{\prime }$ en un point, est tangent \`{a} $\mathcal{F}^{\prime }
$ .
\end{proposition}

\begin{proof}
\bigskip

1) Commen\c{c}ons par montrer que tout champ global feuillet\'{e} transverse
d'un $G-$ feuilletage de Lie minimal ( $\mathcal{M}$, $\mathcal{F}$)d'une
vari\'{e}t\'{e} compacte connexe nul en un point est identiquement nul. \ En
effet \ , dans ces conditions, on sait que l'alg\`{e}bre de Lie structurale $%
\ell $( $\mathcal{M}$, $\mathcal{F}$) de $\mathcal{F}$ et l'alg\`{e}bre de
Lie de $G$ sont isomorphes de dimension la codimension de $\mathcal{F}$.
Ensuite, le feuilletage de Lie $\mathcal{F}$ \'{e}tant \`{a} feuilles
denses, si ($Y_{1},..,Y_{q})$ est un parall\'{e}lisme de Lie transverse de $%
\mathcal{F}$, et si $\ $on consid\`{e}re une m\'{e}trique quasi-fibr\'{e}e
de $\mathcal{F}$ , alors

- l'application " \'{e}valuation en $x$ ", not\'{e}e $ev_{x}$ de $\ell $( $%
\mathcal{M}$, $\mathcal{F}$) dans $\nu _{x}($ $\mathcal{F)\cong }(T_{x}%
\mathcal{F)}^{\bot }$ , d\'{e}finie par $ev_{x}(X)=X(x)$, r\'{e}alise un
isomorphisme canonique d'espaces vectoriels entre $\ell $( $\mathcal{M}$, $%
\mathcal{F}$) et $(T_{x}\mathcal{F)}^{\bot },$

- pour tout $Z\in \ell $( $\mathcal{M}$, $\mathcal{F}$), $Z$ s'\'{e}crit :%
\begin{equation*}
Z=\underset{i=1}{\overset{q}{\sum \xi ^{i}}}Y_{i}
\end{equation*}

o\`{u} les fonctions $\mathcal{F-}$basiques $\xi ^{i}$ sont en fait des
constantes r\'{e}elles. Ceci dit, si $Z$ s'annule en un point donn\'{e}$,$
les $\xi ^{i}$ sont tous nuls et par cons\'{e}quent $Z=0.$

2) Consid\'{e}rons \ maintenant un champ $Z\in \ell $( $\mathcal{M}$, $%
\mathcal{F}$) tel que pour un point $x_{0}$ donn\'{e}, $Z(x_{0})\in T_{x_{0}}%
\mathcal{F}^{\prime }.$ Suivant la d\'{e}composition%
\begin{equation*}
(T\mathcal{F)}^{\bot }=(T\mathcal{F)}^{\bot }\cap T\mathcal{F}^{\prime
}\oplus (T\mathcal{F)}^{\bot }\cap (T\mathcal{F}^{\prime }\mathcal{)}^{\bot }
\end{equation*}

$Z$ s'\'{e}crit $Z=Z_{1}+Z_{2},$ o\`{u} $Z_{1}$et $Z_{2}$ sont des sections
\ globales respectives des sous-fibr\'{e}s $(T\mathcal{F)}^{\bot }\cap T%
\mathcal{F}^{\prime }$ et $(T\mathcal{F)}^{\bot }\cap (T\mathcal{F}^{\prime }%
\mathcal{)}^{\bot }=(T\mathcal{F}^{\prime }\mathcal{)}^{\bot }$ du fibr\'{e}
orthogonal $(T\mathcal{F)}^{\bot }$de $\mathcal{F}.$ Comme $Z(x_{0})\in
(T_{x_{0}}\mathcal{F)}^{\bot }\cap T_{x_{0}}\mathcal{F}^{\prime }$ , alors $%
Z_{2}(x_{0})=0.$ Ce qui implique d'apr\`{e}s le 1) que $Z_{2}$ est
identiquement nul , donc $Z=Z_{1}$ , $i.e.Z$ est tangent \`{a} $\mathcal{F}%
^{\prime }$.
\end{proof}

Ceci \'{e}tant, le th\'{e}or\`{e}me suivant , qui est le r\'{e}sultat
principal de ce travail, d\'{e}termine et classifie les extensions d'un $G$%
-feuilletage de Lie minimal d'une vari\'{e}t\'{e} compacte et connexe.

\begin{theorem}
\label{théo.}Soit ( $\mathcal{M}$, $\mathcal{F}$) un $G$-feuilletage de Lie
minimal d'une vari\'{e}t\'{e} compacte connexe, d'alg\`{e}bre de Lie $%
\mathcal{G}$.

Alors:

1- Il y a une correspondance biunivoque entre les sous-groupes de Lie
connexes de $G$ et les extensions de $\mathcal{F}.$

2- Une extension de $\mathcal{F}$ est un $\frac{\mathcal{G}}{\mathcal{H}}-$%
feuilletage transversalement riemannien \`{a} fibr\'{e} normal trivial, d%
\'{e}finie par une 1- $forme$ vectorielle.

3-Une extension de $\mathcal{F}$ est transveralement homog\`{e}ne( resp. de
Lie) si et seulement si le sous-groupe de Lie de $G$ correspondant est un
sous-groupe ferm\'{e} ( resp. un sous-groupe normal) dans $G.$
\end{theorem}

\begin{proof}
Etant donn\'{e} un $G-$feuilletage de Lie minimal ( $\mathcal{M}$, $\mathcal{%
F}$) de codimension $q>0$ d'une vari\'{e}t\'{e} compacte connexe, d'alg\`{e}%
bre de Lie $\mathcal{G}$,

1. Soit $\mathcal{F}^{\prime }$ une extension de $\mathcal{F}$ , de
codimension $q^{\prime }.$\ Consid\'{e}rons $\widetilde{\mathcal{H}}$
l'ensemble des champs $\mathcal{F-}$feuillet\'{e}s tangents \`{a} $\mathcal{F%
}^{\prime };\widetilde{\mathcal{H}}$ est visiblement une sous-alg\`{e}bre de
l'alg\`{e}bre de Lie structurale $l(\mathcal{M},\mathcal{F)}$ de $\mathcal{F}
$. Pour d\'{e}terminer la dimension de cette sous-alg\`{e}bre, notons, d'apr%
\`{e}s ce qui pr\'{e}c\`{e}de que

- pour tout $x\in \mathcal{M}$ , et pour tout $X\in $ $\widetilde{\mathcal{H}%
}$ , $ev_{x}(X)\in (T_{x}\mathcal{F)}^{\bot }\cap T_{x}\mathcal{F}^{\prime } 
$,

- \ pour tout $u\in (T_{x}\mathcal{F)}^{\bot }\cap T_{x}\mathcal{F}^{\prime
},$ la proposition \ref{pro.2} pr\'{e}c\'{e}dente permet de voir que le
champ $X_{u}=ev_{x}^{-1}(u)$ est dans $\widetilde{\mathcal{H}}.$

Au total , la restriction \`{a} $\widetilde{\mathcal{H}}$ de $ev_{x}$ est un
isomorphisme d'espaces vectoriels de $\widetilde{\mathcal{H}}$ sur $ev_{x}(%
\widetilde{\mathcal{H}})=(T_{x}\mathcal{F)}^{\bot }\cap T_{x}\mathcal{F}%
^{\prime }$ ; ce qui assure par la formule des dimensions que 
\begin{eqnarray*}
\dim \widetilde{\mathcal{H}} &=&\dim (T_{x}\mathcal{F)}^{\bot }\cap T_{x}%
\mathcal{F}^{\prime } \\
&=&\dim (T_{x}\mathcal{F)}^{\bot }+\dim T_{x}\mathcal{F}^{\prime }-\dim
<(T_{x}\mathcal{F)}^{\bot }\cup T_{x}\mathcal{F}^{\prime }> \\
&=&\dim \mathcal{F}^{\prime }\mathcal{-}\dim \mathcal{F} \\
&\mathcal{=}&co\dim \mathcal{F-}co\dim \mathcal{F}^{\prime }
\end{eqnarray*}%
Soit $\omega $ la 1-$forme$ de F\'{e}dida d\'{e}finissant $\mathcal{F}$, et
soit ,en tout point $x$ de $\mathcal{M}$ , $\varpi _{x}$ l'isomorphisme
canonique d'espaces vectoriels \ rendant commutatif le diagramme de
projections%
\begin{equation*}
\begin{array}{cc}
T_{x}\mathcal{M} & \overset{\omega _{x}}{\rightarrow }\mathcal{G} \\ 
\downarrow & \nearrow _{\varpi _{x}} \\ 
\nu _{x}(\mathcal{F)} & 
\end{array}%
\end{equation*}%
et $\sigma $ l'application de$\ \mathcal{G}$ dans $l(\mathcal{M},\mathcal{F)}
$ dans qui \`{a} tout $\lambda \in \mathcal{G}$ associe $\sigma $($\lambda $%
) d\'{e}finie par $\sigma $($\lambda $)$(x)=(\varpi _{x}\circ ev_{x})^{-1}$($%
\lambda );\sigma $ est un isomorphisme lin\'{e}aire de $\mathcal{G}$ sur $l(%
\mathcal{M},\mathcal{F)}$

On remarquera que pour tous $X,Y\in l(\mathcal{M},\mathcal{F)}$,

1) $\sigma ^{-1}$($X$)=$\omega $($X$) ,

2) $\omega $($X$) est une fonction constante sur $\mathcal{M}$ et

3) $\omega \lbrack X,Y]=$ $[\omega (X),\omega (Y)].$

Il en r\'{e}sulte que $\sigma $ est aussi un isomorphisme d'alg\`{e}bres:
c'est cet isomorphisme canonique qui permet d'identifier $\mathcal{G}=Lie(G)$
et l'alg\`{e}bre de Lie structurale $l(\mathcal{M},\mathcal{F)}$ de $%
\mathcal{F}$. \ Ici pour la clart\'{e} de la d\'{e}monstration nous nous
garderons de faire une telle identification.

Ceci \'{e}tant, consid\'{e}rons le syst\`{e}me diff\'{e}rentiel $\mathcal{P}$
d\'{e}fini sur $\mathcal{M}$ par 
\begin{equation*}
\mathcal{P}(x)=T_{x}\mathcal{F\oplus }ev_{x}(\widetilde{\mathcal{H}})\text{
\ }
\end{equation*}%
Soit $\mathcal{X(P)}$ et $\mathcal{X(F)}$ les $\mathcal{A}^{0}(\mathcal{M)-}$%
modules des champs de vecteurs tangents \newline
respectivement \`{a} \ $\mathcal{P}$ et \`{a} $\mathcal{F}$. On a%
\begin{equation*}
\mathcal{X(P)}=\mathcal{X(F)\oplus (A}^{0}(\mathcal{M)\otimes }\widetilde{%
\mathcal{H}})
\end{equation*}%
Comme les champs de vecteurs de $\widetilde{\mathcal{H}}$ sont feuillet\'{e}%
s pour $\mathcal{F}$, cette d\'{e}composition permet de voir que le module $%
\mathcal{X(P)}$ est stable par le crochet et que par suite $\mathcal{P}$ est
un syst\`{e}me diff\'{e}rentiable compl\`{e}tement int\'{e}grable qui d\'{e}%
finit $\mathcal{F}^{\prime }.$ Ainsi la sous-alg\`{e}bre de Lie $\mathcal{H}%
=\sigma ^{-1}(\widetilde{\mathcal{H}})$ de $\mathcal{G}$ et le sous-groupe
de Lie connexe $H$ dans $G$ correspondant sont d\'{e}finis sans ambig\"{u}t%
\'{e} \`{a} partir de l'extension $\ \mathcal{F}^{\prime }.$

R\'{e}ciproquement la donn\'{e}e d'un sous-groupe de Lie connexe $H$ de $G$
permet \ de d\'{e}finir sur $\mathcal{M}$ un syst\`{e}me diff\'{e}rentiel $%
\mathcal{P}$ \ par 
\begin{equation*}
\mathcal{P}(x)=T_{x}\mathcal{F+}\text{ }ev_{x}(\sigma (Lie(H))\text{ (cette
somme est en fait directe).}
\end{equation*}

et le module de champs correspondant \'{e}tant 
\begin{equation*}
\mathcal{X(P)}=\mathcal{X(F)\oplus (A}^{0}(\mathcal{M)\otimes }\sigma
(Lie(H))
\end{equation*}

et $\sigma $ \'{e}tant un isomorphisme d'alg\`{e}bres, il est facile de voir
que le syst\`{e}me diff\'{e}rentiel $\mathcal{P}$ ainsi d\'{e}fini est compl%
\`{e}tement int\'{e}grable, et le feuilletage $\mathcal{F}^{\prime }$qu'il d%
\'{e}finit est bien s\^{u}r une extension de $\mathcal{F}$ puisque pour tout 
$x\in \mathcal{M},$ on a $T_{x}\mathcal{F\subset P}(x)=T_{x}\mathcal{%
\mathcal{F}}^{\prime }$ .

Et la correspondance biunivoque est ainsi \'{e}tablie.

2. - Soit $\mathcal{F}_{H}$ une extension de $\mathcal{F}$ et $\mathcal{H}$
\ la sous-alg\`{e}bre de Lie de $\mathcal{G}$ correspondant \`{a} $H$, soit (%
$e_{1},...,e_{q})$ une base de l'alg\`{e}bre de Lie $\mathcal{G}$ telle que (%
$e_{q^{\prime }+1},...,e_{q})$ engendre $\mathcal{H}$ et $\omega
=\sum\limits_{i=1}^{q}\omega ^{i}\otimes e_{i}$ la $1-forme$ de F\'{e}dida
de $\mathcal{F}$. Puisque en tout point de $\mathcal{M}$ \ les $1-formes$
scalaires $\omega ^{1},...,\omega ^{q}$ sont lin\'{e}airement ind\'{e}%
pendantes $,$ alors les $1-formes$ $\omega ^{1},...,\omega ^{q^{\prime }}$
sont aussi lin\'{e}airement ind\'{e}pendantes partout et la condition de
Mauer- Cartan (*) assure que le syst\`{e}me diff\'{e}rentiel $\omega
^{1}=...=\omega ^{q^{\prime }}=0$ est compl\`{e}tement int\'{e}grable et d%
\'{e}finit donc un feuilletage $\mathcal{F}^{\prime }$de codimension $%
q^{\prime }$ qui n'est rien d'autre que le $\frac{\mathcal{G}}{\mathcal{H}}%
-feuilletage$ d\'{e}fini \ par $\omega .$

Par ailleurs, encore comme $\sigma ^{-1}$($X$)=$\omega $($X$) si $\ X\in 
\widetilde{\mathcal{H}}=\sigma (\mathcal{H)}$ alors pour tout $X\in 
\widetilde{\mathcal{H}}$, $\omega (X)\in \mathcal{H}$. Ce qui permet alors
de voir que pour tout $k,1\leq k\leq q^{\prime }$, et pour tout $\ X\in 
\widetilde{\mathcal{H}},$ $\omega ^{k}(X)=0,i.e.$ $X$ est tangent \`{a} $%
\mathcal{F}^{\prime }.$ Ce qui montre que pour tout $\ x\in M,$ $T_{x}%
\mathcal{F}_{H}=T_{x}\mathcal{F}\mathcal{\oplus }$ $ev_{x}(\widetilde{%
\mathcal{H}})\subset T_{x}\mathcal{F}^{\prime }$ et \`{a} cause d'\'{e}galit%
\'{e} des dimensions, on a $\ $en fait $T_{x}\mathcal{F}_{H}=T_{x}\mathcal{F}%
^{\prime };$ $\mathcal{F}_{H}$ est bien le $\frac{\mathcal{G}}{\mathcal{H}}%
-feuilletage$ associ\'{e} \`{a} $\mathcal{F}$

- R\'{e}ciproquement si $\mathcal{F}^{\prime }$ est un $\frac{\mathcal{G}}{%
\mathcal{H}}-feuilletage$ d\'{e}fini \ par $\omega $ et si $\mathcal{D}$ est
une d\'{e}veloppante de F\'{e}dida de $\mathcal{F}$ d\'{e}finie sur le rev%
\^{e}tement universel $\widetilde{\mathcal{M}\text{ }}$ de $\mathcal{M}$, et
si \ $\widetilde{\mathcal{F}\text{ }}$et $\widetilde{\mathcal{F}^{\prime }}$
sont les feuilletages relev\'{e}s sur $\widetilde{\mathcal{M}\text{ }}$
respectifs de $\mathcal{F}$ et $\mathcal{F}^{\prime },$ alors l'application $%
\mathcal{D}$ \'{e}tant aussi une d\'{e}veloppante du $\frac{\mathcal{G}}{%
\mathcal{H}}-$feuilletage $\mathcal{F}^{\prime },$ par la prop. \ref{prop.},
on a $\widetilde{\mathcal{F}^{\prime }}=\mathcal{D}^{\ast }\mathcal{F}_{G,%
\mathcal{H}}.$ Comme $\widetilde{\mathcal{F}}=\mathcal{D}^{\ast }\mathcal{F}%
_{G,\{e\}}$ et $\mathcal{F}_{G,\{e\}}\subset \mathcal{F}_{G,\mathcal{H}}$,
alors $\widetilde{\mathcal{F}}\subset \widetilde{\mathcal{F}^{\prime }}$ et $%
\ \mathcal{F}^{\prime }$ est une extension de $\mathcal{F}$.

Montrons maintenant que le feuilletage $\mathcal{F}_{H}$ est riemannien.
Pour cela, comme le probl\`{e}me est local, remarquons que:

(**) une extension $\mathcal{F}^{\prime }$ d'un feuilletage riemannien $%
\mathcal{F}$ est un feuilletage riemannien si et seulement si \ la
projection canonique de la restriction de $\mathcal{F}^{\prime }$ \`{a} tout
ouvert \ \`{a} la fois distingu\'{e} pour $\mathcal{F}$ et $\mathcal{F}%
^{\prime }$ est un feuilletage riemannien de la vari\'{e}t\'{e} quotient
locale de $\mathcal{F}$ .

Soit $U$ un ouvert distingu\'{e} \`{a} la fois pour $\mathcal{F}$ et $%
\mathcal{F}_{H},\pi $ et $\pi _{h}$ les projections sur les vari\'{e}t\'{e}s
quotient locales respectives $V$ et $V_{h}$. Alors puisque $\mathcal{F}_{H}$
est une extension de $\mathcal{F},$ il existe une submersion $\theta _{h}$
de $V$ sur $V_{h}$ rendant commutatif le diagramme de projections suivant:%
\begin{equation*}
\begin{array}{cc}
U & \overset{\pi _{h}}{\rightarrow }V_{h} \\ 
^{\pi }\downarrow & \nearrow _{\theta _{h}} \\ 
V & 
\end{array}%
\end{equation*}%
Ce diagramme montre que la projection sur $V$\ de la restriction de $%
\mathcal{F}_{H}$ \`{a} $U$ est un feuilletage simple $(V,\mathcal{F}_{\theta
_{h}})$ puisque d\'{e}fini par la submersion $\theta _{h}$ . Ensuite le
feuilletage $\mathcal{F}$ \'{e}tant un feuilletage de Lie minimal, on sait
que son faisceau central transverse $\mathcal{C(M}$,$\mathcal{F)}$ est
localement trivial et s'identifie aux germes d\'{e}finis par l'alg\`{e}bre
de Lie structurale $l(\mathcal{M},\mathcal{F)}$de $\mathcal{F}$( plus pr\'{e}%
cisement par son alg\`{e}bre de Lie oppos\'{e}e $l(\mathcal{M},\mathcal{F)}%
^{-}$constitu\'{e}e par les champs de vecteurs invariants \`{a} droite). En
supposant en plus que $U$ est un ouvert trivialisant du faisceau, alors le
\textquotedblright sous-faisceau\textquotedblright\ $\mathcal{C}_{H}\mathcal{%
(}U$,$\mathcal{F)=}U\times \sigma (\mathcal{H)}^{-}$ de $\mathcal{C(M}$,$%
\mathcal{F)}$ correspondant \`{a} la sous-alg\`{e}bre de Lie $\sigma (%
\mathcal{H)}$ de $l(\mathcal{M},\mathcal{F)}$ d\'{e}finit , par ses orbites
dans $U$, le feuilletage $\mathcal{F}_{H},_{U}$ restriction de $\mathcal{F}%
_{H}$ \`{a} $U.$ Ainsi on peut regarder les feuilles de $\mathcal{F}%
_{H},_{U} $ comme les orbites du \textquotedblright
sous-faisceau\textquotedblright\ \ $\mathcal{C}_{H}\mathcal{(}U$,$\mathcal{F)%
}$ , $i.e.$ des orbites d'un faisceau de germes de champs de Killing
transverses. Et le feuilletage $\mathcal{F}_{\theta _{h}}$ \'{e}tant la
projection de ces orbites par $\pi $ sur la vari\'{e}t\'{e} quotient locale $%
V$ est alors un feuilletage riemannien dont une m\'{e}trique quasi-fibr\'{e}%
e est la m\'{e}trique projet\'{e}e de la m\'{e}trique transverse de $%
\mathcal{F}$. Et il suit , d'apr\`{e}s la remarque (**), que le feuilletage $%
\mathcal{F}_{H}$ est riemannien.

Ensuite, l'inclusion $T\mathcal{F\subset }T\mathcal{F}_{H}$ induit un
morphisme canonique de fibr\'{e}s vectoriels $\alpha $ de $\vartheta (%
\mathcal{F)}$ sur $\vartheta (\mathcal{F}_{H}\mathcal{)}$ de sorte que le
fibr\'{e} normal de $\mathcal{F}_{H}$ est le fibr\'{e} quotient du fibr\'{e}
normal de $\mathcal{F}$ par son sous-fibr\'{e} $Ker\alpha $. Comme $%
\vartheta (\mathcal{F)}$ en tant que fibr\'{e} normal d'un feuilletage de
Lie est trivial $(\vartheta (\mathcal{F)}\cong M\times \mathcal{G)}$ et
comme $Ker\alpha $ est aussi trivial $\ (Ker\alpha \cong M\times \mathcal{H)}
$ , alors $\vartheta (\mathcal{F}_{H}\mathcal{)}$ est trivialisable et de
section globale ($\sigma (e_{1}),...,\sigma (e_{q^{\prime }}))$. De fa\c{c}%
on pr\'{e}cise $\vartheta (\mathcal{F}_{H})\cong M\times \frac{\mathcal{G}}{%
\mathcal{H}}\cong M\times \mathcal{H}^{\perp }$ \ et o\`{u} $\ \mathcal{G}=%
\mathcal{H}\oplus \mathcal{H}^{\bot }$ et $\mathcal{H}^{\bot }$ l'orthocompl%
\'{e}ment de $\mathcal{H}$ dans l'espace euclidien $\mathcal{G=}T_{e}G$ muni
de la base orthonorm\'{e}e ($e_{1},...,e_{q}).$

3 -\ Soit $\mathcal{F}_{H}$ une extension d'un feuilletage de Lie minimal $\ 
\mathcal{F}$ d'une vari\'{e}t\'{e} compacte $\mathcal{M}$ , $H$ \'{e}tant le
sous- groupe de Lie connexe associ\'{e} \`{a} cette extension. Soit $%
\mathcal{(D}$, $\mathcal{\rho )}$ un d\'{e}veloppement de F\'{e}dida de $%
\mathcal{F}$ sur \ le rev\^{e}tement universel $\widetilde{\mathcal{M}\text{ 
}}$de $\mathcal{M}$. Puisque le feuilletage $\mathcal{F}_{H}$ est aussi un $%
\frac{\mathcal{G}}{\mathcal{H}}-$feuilletage, notons $\widetilde{\mathcal{F}%
_{H}}$ le rel\`{e}vement de $\mathcal{F}_{H}$ sur $\widetilde{\mathcal{M}}$
; $\mathcal{F}_{H}$ est alors de m\^{e}me d\'{e}veloppante que le
feuilletage de Lie $\mathcal{F}$ .

a) En gardant les notations de la prop. \ref{prop.}, il est clair que si le
sous- groupe $H$ est ferm\'{e}, le feuilletage $\mathcal{F}_{G,H}$ $\ i.e.$
le feuilletage de $G$ par les translat\'{e}s de $H$ , est d\'{e}fini par la
submersion canonique $\theta :G\longrightarrow \frac{G}{H},$ et comme d'apr%
\`{e}s la prop. \ref{prop.},$\widetilde{\text{ }\mathcal{F}_{H}}$=$\mathcal{D%
}^{\ast }\mathcal{F}_{G,H},$ il vient que $\mathcal{D}_{H}$=$\mathcal{\theta
\circ D}$ : $\widetilde{\mathcal{M}}$ $\rightarrow \frac{G}{H}$ est une
submersion d\'{e}finissant $\widetilde{\mathcal{F}_{H}}$, \'{e}quivariante
pour la repr\'{e}sentation \ $\rho :\pi _{1}(\widetilde{M})\longrightarrow G$
associ\'{e}e au feuilletage de Lie $\mathcal{F}$; par suite cette extension
est un feuilletage transversalement homog\`{e}ne.

b) R\'{e}ciproquement si $\mathcal{F}_{H}$ est une extension de $\mathcal{F}$
transversalement homog\`{e}ne de vari\'{e}t\'{e} transverse une vari\'{e}t%
\'{e} homog\`{e}ne $T$, comme $\widetilde{\mathcal{F}_{H}}$ est une
extension de$\widetilde{\text{ }\mathcal{F}}$ et que ces deux feuilletages
sont simples( \cite{BLU}, \cite{FED}), alors il existe une submersion $%
\theta $ de $G$ sur $T$ \ telle que la submersion $\mathcal{D}^{\prime }%
\mathcal{=\theta \circ D}$ d\'{e}finit $\widetilde{\mathcal{F}_{H}}$ . Si $%
\mathcal{F}_{\theta }$ est le feuilletage simple donn\'{e}e par la
submersion $\theta ,$ il est clair que%
\begin{equation*}
\mathcal{D}^{\ast }\mathcal{F}_{G,H}=\widetilde{\mathcal{F}_{H}}=\mathcal{D}%
^{\ast }\mathcal{F}_{\theta }
\end{equation*}%
Puisque $\mathcal{D}$ est une submersion surjective on a \'{e}videmment 
\begin{equation*}
\mathcal{F}_{G,H}=\mathcal{DD}^{\ast }\mathcal{F}_{G,H}=\mathcal{DD}^{\ast }%
\mathcal{F}_{\theta }=\mathcal{F}_{\theta },
\end{equation*}%
et il en r\'{e}sulte que le sous-groupe $H$ est la composante connexe de la
fibre $\theta ^{-1}(\theta (e))$ qui contient l'\'{e}l\'{e}ment neutre $e$
de $G.$ Comme cette fibre est un ferm\'{e}e dans $G$, alors $H$ est aussi
une partie ferm\'{e}e de $G.$

Ce qui pr\'{e}c\`{e}de en a) assure que cette extension $\mathcal{F}_{H}$
est un $\frac{G}{H}-$feuilletage transversalement homog\`{e}ne dont $\rho $
serait une r\'{e}pr\'{e}sentation.

En plus, si le sous-groupe $H$ est normal dans $G$, alors le feuilletage $%
\mathcal{F}_{G,H}$ , i.e. le feuilletage de $G$ \ par les translat\'{e}s de $%
H$ est un $\frac{\mathcal{G}}{\mathcal{H}}-$feuilletage de Lie. Comme le
groupe de Lie $G$ est pris connexe et simplement connexe, alors ce
feuilletage est necessairement un feuilletage simple \cite{FED} ; parsuite
le sous-groupe $H$ \'{e}tant la feuille passant par l'\'{e}l\'{e}ment neutre
est ferm\'{e}e et le feuilletage $\mathcal{F}_{G,H}$ d\'{e}fini par la
projection canonique $\theta $ : $G\longrightarrow \frac{G}{H}$ .\ Ensuite
puisque $\theta $ est un mophisme de groupes, alors le feuilletage extension 
$\ \mathcal{F}_{H}$ de $\mathcal{F}$ est un feuilletage de Lie car son
feuilletage relev\'{e} $\widetilde{\mathcal{F}_{H}}$ sur $\widetilde{%
\mathcal{M}}$ \ est d\'{e}fini par la submersion $\mathcal{\theta \circ D}$ 
\'{e}quivariante pour la \ r\'{e}pr\'{e}sentation $\rho ^{\prime }=\theta
\circ \rho :\pi _{1}(\mathcal{M})\longrightarrow \frac{G}{H}.$

R\'{e}ciproquement \ si $\ $l'extension $\mathcal{F}_{H}$ est un feuilletage
de Lie et si \newline
$\mathcal{D}$ : $\widetilde{\mathcal{M}}$ $\rightarrow G$ et $\mathcal{D}%
_{H} $ :$\widetilde{\mathcal{M}}$ $\rightarrow G^{\prime }$ sont des d\'{e}%
veloppantes de F\'{e}dida respectives pour $\mathcal{F}$ et $\mathcal{F}_{H}$%
, alors , on montre comme dans \cite{DIA}, que si $\Gamma $ est le groupe
d'holonomie du feuilletage de Lie $\mathcal{F}$, il existe une submersion $%
\theta $ de $G$ sur $G^{\prime }$ telle que

1)%
\begin{equation*}
\mathcal{D}_{H}\text{ }\mathcal{=\theta \circ D}
\end{equation*}

2) Pour tout $\gamma \in \Gamma $ $,$ le diagramme suivant est commutatif%
\begin{equation*}
\begin{array}{ccccc}
& G & \overset{\theta }{\rightarrow } & G &  \\ 
\gamma & \downarrow &  & \downarrow & \theta (\gamma ) \\ 
& G & \overset{\theta }{\rightarrow } & G & 
\end{array}%
\end{equation*}

$i.e.$ pour tous $\gamma \in \Gamma $ et $g\in G$, 
\begin{equation*}
\theta (\gamma \cdot g)=\theta (\gamma )\cdot \theta (g)
\end{equation*}%
Comme $\Gamma $ est dense dans $G$ et que la restriction de $\theta $ \`{a} $%
\Gamma $ est un homomorphisme de groupes \ qui est continu, alors par
continuit\'{e} , $\theta $ est \'{e}videmment un morphisme de groupes. Dans
ces conditions, $\mathcal{F}_{H}$ \ est aussi une extension de Lie de $%
\mathcal{F}$ correspondant au sous-groupe $K$ composante connexe de l'\'{e}l%
\'{e}ment neutre du sous-groupe normal $\ Ker\theta $ de $G.$ En raison de
la correspondance biunivoque entre extensions et sous-groupes de Lie
connexes , on a necessairement $\ H=K$ et le sous-groupe $H$ est alors
normal .
\end{proof}

\begin{remark}
\label{Rem2}\bigskip
\end{remark}

1- Il d\'{e}coule de cette derni\`{e}re caract\'{e}risation que\textit{\
toute extension de Lie d'un feuilletage homog\`{e}ne minimal d'une vari\'{e}t%
\'{e} compacte \cite{Ghy} est \'{e}galement un feuilletage homog\`{e}ne.}

2- L'hypoth\`{e}se que \textit{le feuilletage de Lie }$\mathcal{F}$\textit{\
est minimal est essentielle} comme le montre l'exemple \ref{Exemple}.

En gardant les notations de la preuve du th\'{e}or\`{e}me pr\'{e}c\'{e}dent,
consid\'{e}rons pour un $G-$feuilletage de Lie minimal $(\mathcal{M},%
\mathcal{F)}$ d'une vari\'{e}t\'{e} compacte et connexe, \ la $1-forme$ $%
\varpi _{H}$ sur $\mathcal{M}$ \`{a} valeurs dans $\mathcal{G}$ d\'{e}finie
par $\varpi _{H}=\alpha \circ \omega ,$ o\`{u} $\alpha $ est la projection
canonique de $\mathcal{G}=\mathcal{H}\oplus \mathcal{H}^{\bot }$ sur $%
\mathcal{H}^{\bot }$ . On v\'{e}rifie que $\varpi _{H}$ est une \'{e}quation
de $\mathcal{F}_{H}$. On peut associer \`{a} cette une $1-forme$ une $2-$ $%
forme$ de courbure $\Omega _{H}=d\omega _{H}+\frac{1}{2}[\omega _{H},\omega
_{H}]$ d\'{e}finie sur $\mathcal{M}$ \`{a} valeurs dans $\mathcal{G}$, o\`{u}
$\omega _{H}=j\circ \varpi _{H}$, $j$ \'{e}tant l'injection canonique de $%
\mathcal{H}^{\bot }$ dans $\mathcal{G}$.

En partant de la formule classique,%
\begin{equation*}
\Omega _{H}(X,Y)=X\omega _{H}(Y)-Y\omega _{H}(X)-\omega _{H}[X,Y]+[\omega
_{H}(X),\omega _{H}(Y)],
\end{equation*}%
et en remarquant que, pour tout $\ X\in l(\mathcal{M},\mathcal{F)}$ , $%
\omega _{H}(X)$ est une fonction constante$,$ alors un calcul facile permet
de voir que:

1) $\Omega _{H}(X,Y)=0$ si $X$ et $Y$ sont tangents \`{a} \ $\mathcal{F}%
^{\prime },$

2) $\Omega _{H}(X,Y)=-\alpha \sigma ^{-1}[X,Y]$ si $X$ $\in \sigma (\mathcal{%
H})=$ $\widetilde{\mathcal{H}}$ et $Y\in \sigma (\mathcal{H}^{\bot })$

3) $\Omega _{H}(X,Y)=(1-\alpha )\sigma ^{-1}[X,Y]$ si $\ X\in \sigma (%
\mathcal{H}^{\bot })$ et $Y\in \sigma (\mathcal{H}^{\bot })$

La $2-forme$ de courbure $\Omega _{H}$ \'{e}tant une application $\mathcal{A}%
^{0}(\mathcal{M})$- bilin\'{e}aire de $\mathcal{X(M)\times X(M)}$ dans $%
\mathcal{G},$ les relations pr\'{e}c\'{e}dentes d\'{e}terminent parfaitement
la\textit{\ }$2-forme$ $\Omega _{H}.$

En partant de cette $1-forme$ $\omega _{H},$ on obtient les caract\'{e}%
risations suivantes d'une extension de Lie d'un feuilletage de Lie minimal
d'une vari\'{e}t\'{e} compacte.

\begin{corollary}
Si $\mathcal{F}_{H}$ est une extension d'un G-feuilletage de Lie minimal
d'une vari\'{e}t\'{e} compacte connexe, les assertions suivantes sont \'{e}%
quivalentes.

1. $\mathcal{F}_{H}$ est une extension de Lie

2. $H$ est un sous-groupe normal

3. $\omega _{H}$ est basique pour $\mathcal{F}_{H}$

4. $\Omega _{H}$ = $d\omega _{H}+\frac{1}{2}[\omega _{H},\omega _{H}]$ est
basique pour $\mathcal{F}_{H}$

5. $\Omega _{H}=0.$
\end{corollary}

On notera que lorsque l'extension $\mathcal{F}_{H}$ est de Lie , la $1-$ $%
forme$ $\omega _{H}$ explicit\'{e}e ci-dessus est sa $1-$ $forme$ de F\'{e}%
dida \cite{FED}.

\begin{corollary}
\label{Corol.}Soit ( $\mathcal{M}$, $\mathcal{F}$) un G-feuilletage de Lie
minimal d'une vari\'{e}t\'{e} compacte connexe.

Si le groupe fondamental $\pi _{1}(\mathcal{M})$ \ est virtuellement r\'{e}%
soluble ( resp. ab\'{e}lien) toute extension de $\mathcal{F}$ est
transversalement homog\`{e}ne(resp. de Lie).

En particulier:

1) toute extension d'un feuilletage lin\'{e}aire minimal du tore est \'{e}%
galement lin\'{e}aire,

2) toute extension d'un flot riemannien minimal d'une vari\'{e}t\'{e}
compacte est conjugu\'{e}e \`{a} un feuilletage lin\'{e}aire du tore.
\end{corollary}

\begin{proof}
En consid\'{e}rant sur$\mathcal{M}$ une m\'{e}trique quasi-fibr\'{e}e
quelconque pour $\mathcal{F}$, cette m\'{e}trique \'{e}tant compl\`{e}te et $%
\pi _{1}(\mathcal{M})$ virtuellement r\'{e}soluble, et comme par ailleurs
toute extension de $\mathcal{F}$ est un feuilletage riemannien minimal (th%
\'{e}o.\ref{théo.}) alors d'apr\`{e}s le th\'{e}o.\ref{Hae} cette extension
est aussi un feuilletage transversalement homog\`{e}ne.

Si en plus $\pi _{1}(\mathcal{M})$ est ab\'{e}lien , alors le groupe
d'holonomie de $\mathcal{F}$ \'{e}tant ab\'{e}lien et dense dans $G$ , il
vient par continuit\'{e} de la loi de groupe, que le groupe de Lie $G$ est ab%
\'{e}lien et par suite toute extension de $\mathcal{F}$ est de Lie.

Pour le reste cela tient de la partie 1) de la remarque \ref{Rem2} et du r%
\'{e}sultat de Yves Carri\`{e}re sur les flots riemanniens minimaux \cite%
{CAR}.
\end{proof}

\ \ Ensuite, en rappelant qu'un feuilletage $\mathcal{F}$ est dit \textit{%
dense} dans un autre feuilletage $\mathcal{F}^{\prime }$ si toute feuille de 
$\mathcal{F}$ est dense dans la feuille de $\mathcal{F}^{\prime }$ qui la
contient, on a :

\begin{corollary}
1) Si une vari\'{e}t\'{e} compacte \`{a} groupe fondamental virtuellement r%
\'{e}soluble supporte un $G-$feuilletage de Lie minimal, alors tout
sous-groupe de Lie connexe de $G$ est ferm\'{e}; et toute extension de ce
feuilletage est transversalement homog\`{e}ne.

2) Si un flot $\mathcal{\ }$riemannien ( $\mathcal{M}$, $\mathcal{F}$) d'une
vari\'{e}t\'{e} compacte est dense dans une de ses extensions $\mathcal{F}$'
alors la restriction de $\mathcal{F}$' \`{a} l'adh\'{e}rence de l'une
quelconque de ses feuilles et conjugu\'{e}e \`{a} un feuilletage lin\'{e}%
aire.
\end{corollary}

\begin{proof}
\bigskip 1) En effet au sous-groupe de Lie $H$ de $G$, il correspond d'apr%
\`{e}s \ le th\'{e}o.\ref{théo.} une extension de ce $G-$feuilletage de Lie
minimal qui par le th\'{e}or\'{e}me Haefliger sus-cit\'{e} est un
feuilletage transversalement homog\`{e}ne ; \ ce qui , d'apr\`{e}s le th\'{e}%
o.\ref{théo.} encore implique que $H$ \ est ferm\'{e}.

2) Ceci r\'{e}sulte imm\'{e}diatement de \cite{CAR} et du corollaire \ref%
{Corol.} pr\'{e}c\'{e}dent.
\end{proof}

Par analogie au th\'{e}or\`{e}me de Molino sur les feuilletages
transversalement parall\'{e}lisables \cite{MOL}, le th\'{e}o \ref{théo.}
permet de dire en

\begin{conclusion}
Si $\mathcal{F}^{\prime }$ est une extension d'un feuilletage
transversalement \newline
parall\'{e}lisable $\mathcal{F}$ d'une vari\'{e}t\'{e} compacte connexe,\
telle que $\mathcal{F}$ est dense dans $\mathcal{F}^{\prime }$ alors:

1 - les adh\'{e}rences des feuilles de $\mathcal{F}^{\prime }$ $forme$nt une
fibration localement triviale \'{e}gale \`{a} la fibration basique de $%
\mathcal{F}$,

2 - les feuilles de $\mathcal{F}^{\prime }$ sont les orbites \ d'un
sous-faisceau du faisceau transverse central de $\mathcal{F}$ , de
fibre-type une sous-alg\`{e}bre de Lie oppos\'{e}e \`{a} une sous-alg\`{e}%
bre de l'alg\`{e}bre de Lie structurale de $\mathcal{F}$,

3 - le feuilletage $\mathcal{F}^{\prime }$ est transversalement riemannien
et \`{a} fibr\'{e} normal trivial,

4 - la restriction de $\mathcal{F}^{\prime }$ \`{a} l'adh\'{e}rence d'une
feuille est un $\frac{\mathcal{G}}{\mathcal{H}}-$feuilletage associ\'{e} au $%
\mathcal{G-}$ feuilletage de Lie d\'{e}fini par $\mathcal{F}$ \ dans cette
restriction.
\end{conclusion}

\newpage

\bibliographystyle{AABBRV}
\bibliography{acompat}

\end{document}